\documentclass[12pt]{article}
\usepackage{amssymb, amscd}
\newcommand{\C}{\mathbb{C}}

\newcommand{\N}{\mathbb{N}}

\newcommand{\Q}{\mathbb{Q}}

\newcommand{\Z}{{\mathbb{Z}}}
\newcommand{\F}{\mathbb{F}}
\newcommand{\lra}{\longrightarrow}
\newcommand{\ra}{\rightarrow}

\newcommand{\whl}{\widehat{L}}
\newcommand{\wha}{\widehat{A}}

\newcommand{\whm}{\widehat{M}}

\newcommand{\whj}{\widehat{J}}

\newcommand{\id}{\mbox{\rm id}}

\begin{document}
\def\flexleftarrow#1{\mathtop{\makebox[#1]{\leftarrowfill}}}
\begin{center}
{\bf Curves of genus \, \boldmath$g$\unboldmath \, on an abelian variety of dimension 
\boldmath$g$\unboldmath}\\[1cm]
H. Lange  and E. Sernesi
\\[2cm]
\end{center}
{\small {\bf Abstract.} In this paper we prove a general theorem concerning the number of translation classes of curves of genus
$g$ belonging to a fixed cohomology class in a polarized abelian variety of dimension $g$. For $g = 2$ we recover results 
of G\"ottsche and Bryan-Leung. For $g = 3$ we deduce explicit numbers for these classes.}\\[1cm]

\noindent
{\bf 1. Introduction.}
\vspace{1ex}
\newline
Let $(A,L)$ be a polarized abelian variety of dimension $g$ and of type $(d_1, \ldots , d_g)$ over the field $\C$ of complex numbers.
Two curves, $C_1$ and $C_2$ in $A$ are called {\it equivalent}, if there is an $x \in A$ such
that $C_1 = t_x C_2$. Here $t_x: A \to A$ denotes the translation by $x$.
We call a class of equivalent curves in $A$ a {\it translation class} of curves. Equivalent curves  in $A$
define the same cohomology class in $H^{2g-2} (A, \Z)$.  Let $r \in \Q, r > 0$ such that the cohomology class $ r \bigwedge^{g-1}c_1(L)$ is
contained in $H^{2g-2} (A, \Z)$. It is easy to see that the number of translation classes of irreducible
reduced curves of genus $g$ in the class $r \bigwedge^{g-1} c_1(L)$ is finite. The problem is to 
compute this number. 
\medskip

\noindent
In the case of an abelian surface this problem is equivalent to computing the number
of curves of genus $2$ in a linear system $|L|$. In fact, if $L$ is of type $(d_1, d_2)$,
then $|L|$ contains exactly $d^2_1 \cdot d^2_2$ curves which are equivalent to a given curve
of genus $2$. The number of curves of genus $2$ in $|L|$ has been computed in the special case
of a simple abelian surface of type $(1,d)$ by G\"{o}ttsche and Debarre in [G] and [D] and finally by Bryan and Leung in
[BrL], applying intersection cohomology methods.
\medskip

\noindent
In this paper we prove a general theorem (Theorem 2.1) expressing the number of equivalence classes of curves
of genus $g$ in the class $r \bigwedge^{g-1} c_1(L) \in H^{2g-2} (A, \Z)$ in terms of maximal
isotropic subgroups of Jacobian type of a finite symplectic space associated to $L$. For the definitions see 
Section 2.
In dimension $g \geq 4$ it is difficult to decide whether a given maximal isotropic subgroup is of Jacobian type.
However any maximal isotropic subgroup of a simple abelian variety of dimension 2 and 3 is of Jacobian type.
For these abelian varieties the theorem can be applied to actually compute the number of translation classes of curves in a fixed 
cohomology class. For abelian surfaces we obtain the results of G\"{o}ttsche, Debarre and Bryan-Leung in a slightly more general form.
The main application of the theorem is dimension 3, where we compute the number of translation classes of curves of genus 3 
for any minimal cohomology class in $H^4 (A, \Z)$ and any simple abelian threefold.  
It would be interesting to find
generating functions for these numbers.\\

\noindent
We would like to thank the referee for suggesting some improvements in the proof of Proposition 4.4.\\
\bigskip

\noindent
{\bf 2. Statements of the results.}
\vspace{1ex}
\newline
First we recall some preliminaries.  By a curve we will always mean a complete reduced irreducible curve. 
By its genus we understand its geometric genus, i.e. the genus of its normalization.
A polarization on an abelian variety $A$ is by definition the first Chern class
of an ample line bundle $L$ on $A$. By abuse of notation we denote the polarization by $L$ instead of $c_1(L)$ 
and consider it as a line bundle.
For any polarization $L$ on $A$ the group $K(L)$ is
defined to be the kernel of the isogeny 
$$
\phi_L : A \rightarrow  \widehat{A}, 
\,\, x \mapsto t^{\ast}_x L
\otimes L^{-1}.
$$ 
If $L$ is of type $(d_1, \ldots , d_g)$, then $K(L) \simeq (\Z / d_1 \Z \times \ldots \times
\Z/ d_g \Z)^2$. The group $K(L)$ admits a nondegenerate (multiplicative) alternating form
$e^L : K(L) \times K(L) \ra \C^{\ast}$, so that one can speak of totally isotropic subgroups of
$K(L)$ which we will shortly call {\it isotropic subgroups}. According to [CAV], Corollary 6.3.5 there is a canonical 
bijection between the sets of

\begin{itemize}
\item[(i)] maximal isotropic subgroups of $K(L)$ and
\item[(ii)] isogenies $f : (A,L) \ra (B,M)$ of polarized abelian varieties (i.e. with $f^{\ast} M$ algebraic equivalent to 
 $L$) onto a principally polarized abelian variety $(B,M)$.
\end{itemize}
We call a maximal isotropic subgroup {\it of Jacobian type} if the associated principally
polarized abelian variety $(B,M)$ is the Jacobian of a smooth projective curve.\\
\noindent
For any polarization $L$ on $A$ there is a unique {\it dual polarization} $\whl$ on the
dual abelian variety $\wha$, characterized by the following 2 equivalent properties (see
[BL]).
\begin{itemize}
\item[(i)] $\phi^{\ast}_L \whl \approx L^{d_1 d_g}$
\item[(ii)] $\phi_{\whl} \phi_L = d_1 d_g \id_A$.
\end{itemize}
Here ``$\approx$'' denotes algebraic equivalence. The polarization $\whl$ is of type  $\left( d_1, \frac{d_1 d_g}{d_{g-1}}, \ldots , \frac{d_1
d_g}{d_2}, d_g \right)$, if $L$ is of type $(d_1, \ldots , d_g)$. Moreover, one has 
biduality:  $\widehat{\whl} \approx L$. The main result of the paper is the following theorem:\\

\noindent
{\bf Theorem 2.1.} {\it Let $(A,L)$ be a polarized abelian variety of dimension $g \geq 2$ and type
$(d_1, \ldots , d_g)$ and $r \in \Q, r > 0$, such that $r \bigwedge^{g-1} c_1(L) \in H^{2g-2} (A,
\Z)$. There is a canonical bijection between the sets of}
\begin{itemize}
\item[(1)] {\it translation classes of curves of genus $g$ in the class $r \bigwedge^{g-1} c_1(L)$.}
\item[(2)] {\it maximal isotropic subgroups in $K(\whl^{r(g-1)! d_2 \cdots d_{g-1}})$ of 
Jacobian
type.}
\end{itemize}
{\it Here for $g=2$ the empty product $d_2 \cdots d_{g-1}$ is considered to be equal to 1}.\\

Note that, since a general abelian variety of dimension $\geq 4$ is not isogenous to a Jacobian, for such an $(A,L)$ 
the cardinalities of the sets (1) and (2) of the theorem are both equal to 0. In particular, a generic 
polarized abelian variety of dimension $g \geq 4$ does not contain any curve of genus $g$.
On the other hand it is difficult to decide whether a given polarization in dimension $\geq 4$ is of Jacobian type or not.
Consequently the main interest of the theorem is for abelian varieties of dimensions 2 and 3.\\  

\medskip
\noindent
It is well known (see e.g. [CAV], Corollary 11.8.2) that any simple principally polarized abelian variety of dimension
2 or 3 is the Jacobian of a nonsingular curve. Hence for $g = 2$ and 3 any maximal isotropic subgroup of a simple 
abelian variety is of Jacobian type.
So in these cases we can be more precise.\\

\noindent
{\bf Corollary  2.2:} {\it Let $(A,L)$ be a simple abelian surface of type $(d_1, d_2)$. The number
of curves of genus $2$ in the linear system $|L|$ is $d^2_1 d^2_2 \cdot \# \{\mbox{maximal
isotropic subgroup of}\,\,  K(L)\}$.}
\smallskip

\noindent 
Note that if $L$ is of type $(1,d)$, the number of maximal isotropic subgroups is
$\sigma (d) = \sum_{m|d} m$ (see Proposition 4.3).
\medskip

\noindent
{\it Proof.} It remains to show that the number of curves in $|L|$ equivalent to a given curve
$C \in |L|$ of genus 2 is $d^2_1 d^2_2$. But if $C'$ is equivalent to $C$ then $C' =
t^{\ast}_x C$ for some $x \in A$. This implies $x \in K(L)$. Since $K(L)$ is of order $d^2_1
d^2_2$, it suffices to show that $t^{\ast}_x C \not= C$ for any $x \in A, x \not= 0$. But if
$t^{\ast}_x C =C$, then $t^{\ast}_x$ has to map the finite set of singularities of $C$ into
itself. This is impossible for $x \not= 0$, if we replace $C$ by a suitable
translate. In fact, replacing $C$ by $t^*_{x_0}C$ for a general $x_0$, we obtain the identity 
$t^*_{x-x_0}(t^*_{x_0}C) = t^*_{x_0}C$ with $x-x_0$ not a division point and such a translation $t_{x -x_0}$
cannot stabilize a finite set of points. \hfill $\square$
\medskip

\noindent
Denote by $K(d_1, \ldots , d_g)$ the group $(\Z / d_1 \Z \times \ldots \times \Z / d_g \Z)^2$
together with the standard symplectic form $e: K(d_1, \ldots , d_g)^2 \ra \C^{\ast}$ defined
as follows. If $f_1, \ldots , f_{2g}$ denotes the standard generators of $K(d_1, \ldots ,
d_g)$, then
$$
e(f_{\mu}, f_{\nu}) = 
\left\{
\begin{array}{ccc}
\exp \left( - \frac{2 \pi i}{d_{\nu}}\right) & \mbox{if}& \mu = g+ \nu\\
\exp \left( \frac{2 \pi i}{d_{\nu}} \right)  & \mbox{if} & \nu = g+ \mu\\
1 && \mbox{otherwise.}
\end{array}
\right.
$$
{\bf Corollary 2.3:} {\it Let $(A,L)$ be a simple abelian threefold of type $(d_1, d_2, d_3)$. The
number of translation classes of curves of genus $3$ in the class $r
\bigwedge^2 c_1(L)$ for some $r \in
\Q, r > 0$ is equal to the number of maximal isotropic subgroups of $K(2r d_1 d_2, 2 r d_1
d_3, 2r d_2 d_3$)}.
\bigskip

\noindent
In particular the number of equivalence classes of curves of genus $3$ in $\frac{1}{2 d_1 d_2}
\bigwedge^2 c_1(L)$ is equal to the number of maximal isotropic subgroups of $K \left( \frac{d_3}{d_2},
\frac{d_3}{d_1} \right)$. If $L$ is of type $(d, d, d)$, then there is exactly one class of
curves of genus 3 in the class $\frac{1}{2d^2} \bigwedge^2 c_1(L)$.
\vspace{1cm}
\newline
{\bf 3. Proof of the theorem.} 
\vspace{1ex}
\newline
We need some preliminaries. Let $(A,L)$ be a polarized abelian
variety of dimension $g$ and $C$ a curve of genus $g$ in $A$ generating $A$. Then we have the
following diagram
$$
\begin{array}{ccc}
\widetilde{C} & \hookrightarrow & J = J(\widetilde{C})\\
&&\\
{\scriptstyle \nu} \downarrow && \downarrow {\scriptstyle f}\\
&&\\
C & \hookrightarrow & A
\end{array}
\eqno(1)
$$
where $\nu$ denotes the normalization map, $J$ the Jacobian of $\widetilde{C}$, and $f$ the
homomorphism induced by the composed map $\widetilde{C}
\stackrel{\nu}{\rightarrow} C \hookrightarrow A$. Note that $f$ is an isogeny, since $\widetilde{C}$ is
of genus $g$ and $C$ generates $A$ as an abelian variety. Let $M$ denote a line bundle on $J$
defining the canonical principal polarization on $J$.
\medskip

\noindent
{\bf Lemma 3.1:} {\it For any $e \in \N$ the following statements are equivalent:}
\begin{itemize}
\item[(i)] \quad $f^{\ast}L \equiv M^{e \deg L}$.
\item[(ii)] \quad $[C] = \frac{e}{(g-1)!} \bigwedge^{g-1} c_1(L) \quad \mbox{\it in} \quad
H^{2g-2} (A, \Z)$.
\end{itemize}
Here $\equiv$ denotes numerical equivalence and the {\it degree of} $L$ is defined as $\deg L
= d_1 \ldots d_g$ if $L$ is of type $(d_1, \ldots , d_g)$.
\medskip

\noindent
{\it Proof.} Note first that (i) is equivalent to $\hat{f} \phi_L f = e \deg L \cdot 1_J$, 
i.e. to the commutativity of the diagram
$$
\begin{array}{ccc}
A& \stackrel{\phi_L}{\longrightarrow }& \wha\\
&&\\
{\scriptstyle f}\uparrow && \downarrow {\scriptstyle \hat{f}}\\
&&\\
J & \stackrel{(e \deg L)_J}{\longrightarrow}  & J,
\end{array}
\eqno(2)
$$
where we identify $J = \widehat{J}$ under the isomorphism $\phi_M : J \ra \widehat{J}$.
We claim that (2) is equivalent to
$$
f \hat{f} \phi_L = e \deg L \cdot 1_J. \eqno(3)
$$
For the proof note that being an isogeny $f: J \to A$ is an isomorphim in $Hom_{\Q}(J,A) = Hom(J,A) \otimes_{\Z} \Q$.
Hence (2) is equivalent to
$$
\hat{f} \phi_L = e \deg L \cdot f^{-1}. 
$$
which in turn is equivalent to (3).
\medskip

\noindent 
Recall from [CAV] Section 5.4 that for any cycles $V$ and $W$ of complementary dimensions one defines
an endomorphism $\delta (V,W)$ of $A$ by
$$
\delta (V,W)(x) = S(V \cdot (t^{\ast}_x W - W))
$$
where $S (V \cdot (t^{\ast}_x W - W)) = \sum^n_{i=1} r_i x_i \in A$ if $V \cdot (t^{\ast}_x W - W) = \sum^n_{i=1} r_i x_i $ as a
zero-cycle. Moreover $\delta (V,W)$ depends only on the algebraic equivalence 
classes of $V$ and $W$. In particular $\delta(C, L)$ is a well defined endomorphism of $A$.\\
\noindent
Now assume (3) holds. Then
\begin{eqnarray*}
\delta (C,L) &=& - f \hat{f} \phi_L \quad \hspace*{7.5cm}\mbox{by [CAV], Prop. 11.6.1}\\
&=& -e \deg L \cdot 1_A \quad \hspace*{9.5cm}\mbox{by (3)}\\
&=& \delta \left( \frac{e}{(g-1)!} \bigwedge^{g-1} c_1(L), L) \right) \quad \hspace*{4cm}
\mbox{by [CAV], Prop. 5.4.7}
\end{eqnarray*}
Since $L$ is an ample line bundle this implies using [CAV], Theorem 11.6.4 that $C$ is numerically equivalent to
$\frac{e}{(g-1)!} \bigwedge^{g-1} c_1(L)$   which is equivalent to (ii),
since numerical and homological equivalence coincide on an abelian variety.\\
\noindent
Conversely (ii)
implies
$$
\delta (C,L) = \frac{e}{(g-1)!} \delta (\bigwedge^{g-1} c_1(L),L)
$$
according to a theorem of Matsusaka (see [M]). But as above $\delta (C,L) = - f \hat{f}
\phi_L$ and $\delta (\bigwedge^{g-1} c_1(L),L) = - (g-1)! \deg L$. So this implies (3), which completes
the proof of Lemma 3.1. \hfill $\square$
\bigskip

\noindent
{\bf Lemma 3.2:} {\it Let $(A,L)$ be a polarized abelian variety of type $(d_1, \ldots , d_g)$ and
$(B,M)$ a principally polarized abelian variety. If $g : (A,L) \ra (B,M)$ is an isogeny of
polarized abelian varieties, i.e. $g^{\ast} M \approx L$, then} $\hat{g}^{\ast} \whl \approx
\whm^{d_1d_g}$.
\medskip

\noindent
{\it Proof.} The polarization $M$ being principal, we may identify $B = \widehat{B}$ such that
$\phi_M = 1_B$. Then $g^{\ast}M \approx L$ implies $\phi_{g^{\ast}M} = \phi_L$ and the following
diagram commutes
$$
\begin{array}{ccccc}
A & \stackrel{\phi_L}{\longrightarrow} & \widehat{A}& \stackrel{\phi_{\whl}}{\longrightarrow} & A\\
&&&&\\
{\scriptstyle g} \downarrow &&  \uparrow {\scriptstyle{\hat{g}}} && \downarrow
{\scriptstyle g}\\
&&&&\\
B &\stackrel{1_B}{\longrightarrow}&B& \stackrel{\phi_{\hat{g}^\ast \whl}}{\longrightarrow} & B
.
\end{array}
$$
By definition of $\whl$ we have
$$
\phi_{\whl} \circ \phi_L = d_1 d_g 1_A.
$$
Hence 
$$
\phi_{\hat{g} \ast \whl} \circ g = d_1 d_g g .
$$
Since $g$ is an isogeny, this gives
$$
\phi_{\hat{g}^{\ast} \whl} = d_1 d_g \cdot 1_B.
$$
On the other hand 
$$
\phi_{\whm^{d_1 d_g}} = d_1 d_g \phi_{\whm} = d_1 d_g \cdot 1_B.
$$
Hence $\phi_{\hat{g} \ast \whl} = \phi_{\whm^{d_1 d_g}}$ which implies the assertion. \hfill
$\square$
\bigskip

\noindent
{\bf Proof of Theorem 2.1:} Let $C$ be a curve of genus $g$ in the class $r
\bigwedge^{g-1} c_1(L)$ in $H^{2g-2} (A, \Z)$. Denote for abbreviation $\beta = r (g-1)! \deg L$.
According to Lemma 3.1
$$
f^{\ast} L \approx M^{\beta}
$$
and the following diagram commutes
$$
\begin{array}{ccc}
A & \stackrel{\phi_L}{\longrightarrow} & \wha\\
{\scriptstyle f} \uparrow && \downarrow {\scriptstyle \hat{f}}\\
J & \stackrel{\beta_J}{\longrightarrow} & J.
\end{array}
$$
Hence 
$$
f^{\ast} \phi_L^{\ast} \hat{f}^{\ast}M = \beta^{\ast}  M \approx M^{\beta^2} \approx f^{\ast}
L^{\beta}.
$$
This implies
$$
\phi^{\ast}_L \hat{f}^{\ast} M \approx L^{\beta}.
$$
On the other hand by definition of $\whl$ we have $\phi^{\ast}_L \whl \approx L^{d_1 d_g}$ and
thus
$$
\phi^{\ast}_L \left( \whl^{\frac{\beta}{d_1 d_g}} \right) \approx L^{\beta}.
$$
It follows
$$
\phi^{\ast}_L \hat{f}^{\ast} M \approx \phi^{\ast}_L \whl^{\frac{\beta}{d_1 d_g}}.
$$
This implies $\hat{f}^{\ast} M \approx \whl^{r (g-1)! d_2 \cdots d_{g-1}}$.
Hence the isogeny $\hat{f}$ corresponds to a maximal isotropic subgroup in $K\left(
\whl^{r(g-1)!d_2 \cdots d_g}\right)$ of Jacobian type. It is clear that equivalent curves
lead to the same maximal isotropic subgroup.\\
\noindent
Conversely let $K$ be a maximal isotropic subgroup of $K(\whl^{r (g-1)!d_2 \cdots d_{g-1}})$
of Jacobian type. Then there is an isogeny 
$$
\hat{f} : \wha \rightarrow \whj = \wha/K
$$
such that
$$
\hat{f}^{\ast} \whm \approx \whl^{r(g-1)!d_2 \cdots d_{g-1}}.
$$
Let $M$ denote the dual polarization of $\widehat{M}$ on $J = \widehat{\whj}$. $(J,M) \simeq (\whj, \whm)$ is
the Jacobian of a smooth curve $\widehat{C}$ of genus $g$. Let $f : J \ra A$ denote the dual
isogeny of $\hat{f}$. Then $C := f(\widetilde{C})$ is a curve of genus $g$ on $A$ and $f |
\widetilde{C} : \widetilde{C} \ra C$ is birational, since $C$ generates the abelian variety
$A$.\\
\noindent
We have to show that $C$ is in the class $r \bigwedge^{g-1} c_1(L) \in H^{2g-2} (A, \Z)$.
Note first that $\whl^{r(g-1)!d_2 \cdots d_{g-1}}$ is of type
$$
r(g-1)!\left( \frac{\deg L}{d_g} , \frac{\deg L}{d_{g-1}}, \ldots , \frac{\deg L}{d_1}
\right).
$$
In particular we have
$$
\deg f = r^g ((g-1)!)^g (\deg L)^{g-1}.
$$
So Lemma 3.2 applied to $g = \hat{f}$ yields
$$
f^{\ast} L^{r(g-1)!d_2 \cdots d_{g-1}} = \hat{\hat{f}}^{\ast} ((\whl^{r(g-1)!d_2 \cdots
d_{g-1}})^{\wedge}) \approx \whm^{ \frac{((r(g-1)! \deg L)^2}{d_1 d_g}}
$$
and thus
$$
f^{\ast} L \approx \whm^{r(g-1)! \deg L}.
$$
On the other hand $f | \widetilde{C} : \widetilde{C} \ra C$ is birational, which gives
\begin{eqnarray*}
f^{\ast}[C] &=& \sum_{x \in K} t^{\ast}_x [ \widetilde{C}]\\
& \equiv& \deg f \cdot [ \widetilde{C}]\\
&=& r^g ((g-1)!)^{g-1} (\deg L)^{g-1} \cdot \bigwedge^{g-1} c_1(\whm)\\
&=& r \bigwedge^{g-1} (r (g-1)! \deg L \cdot c_1(\whm ))\\
&=& r \bigwedge^{g-1} (f^{\ast} c_1(L))\\
&=& f^{\ast} (r \bigwedge^{g-1} c_1(L)).
\end{eqnarray*}
But for an isogeny $f^{\ast}$ is injective on cohomology, implying 
$$
[C] = r \bigwedge^{g-1} c_1(L).
$$
It is easy to see that both maps are inverse to each other. \hfill $\square$\\
\bigskip

\noindent
{\bf 4. Maximal isotropic subgroups.}
\vspace{1ex}
\newline
\noindent
Suppose $d_1,...,d_g$ are positive integers with $d_i|d_{i+1}$ for $i = 1,...,g-1$ and let $K(d_1,...,d_g) 
= (\Z/d_1 \Z \times ... \times \Z/d_g \Z)^2$ denote the finite group with the standard symplectic 
form $e(\,.\,,\,.\,)$ of Section 2. Recall that a subgroup $H \subset K(d_1,...,d_g)$ is called {\it (totally)
isotropic} if $e(h, h')=1$ for all $h, h' \in H$. An isotropic subgroup of $K(d_1,...,d_g)$ is 
{\it maximal isotropic}
if and only if it is of order $d:= d_1 \ldots d_g$. Let 
$$
\nu (d_1,...,d_g)
$$
denote the number of maximal isotropic 
subgroups of $K(d_1,...,d_g)$. In this section we compute the number $\nu (d_1,...,d_g)$ in some cases. 
First note that from the definition of the symplectic form $e(\, \cdot \, , \, \cdot \, )$ we immediately obtain\\

\noindent
{\bf Proposition 4.1:} {\it Let $d'_1,...,d'_g$ be another set of positive integers with $d'_i|d'_{i+1}$ for
$i=1,...,g$. If $(d_g,d'_g) = 1$ then there is a symplectic isomorphism
$$
K(d_1d'_1, ... ,d_gd'_g) \simeq K(d_1,...,d_g) \times K(d'_1,...,d'_g).
$$
In particular $\nu (d_1d'_1,...,d_gd'_g) = \nu (d_1,...,d_g) \cdot \nu (d'_1,...,d'_g)$}.\\

\noindent
{\bf Proposition 4.2:} {\it $\nu (p,...,p) = \prod^g_{i=1} (p^i +1)$ for any prime number $p$.}
\medskip

\noindent
{\it Proof:} According to Witt's theorem the group $Sp(2g,p)$ acts transitively on the set of
maximal isotropic subgroups of $K(p, \cdots ,p)$. Hence
$$
\nu (p, \cdots ,p) = \frac{|Sp(2g,p)|}{| {\rm Fix} K_0|}
$$
where $K_0$ denotes the isotropic subgroup generated by $f_1, \ldots , f_g$, if $f_1, \ldots
, f_{2g}$ denotes a symplectic basis of $K(p, \cdots ,p)$. An element
$\left( \begin{array}{ll}
\alpha & \beta \\ 
\gamma & \delta
\end{array}
\right)
 \in Sp(2g,p)$ is contained in $K_0$ if and only if
$\alpha =0, \,\,\, \beta \,\,= \,\, - ^t\!\gamma^{-1}$ and $\gamma \,\, ^t\!\delta \,\, = \,\, \delta \,\, ^t\!\gamma$. 
Hence
$$
K_0 \simeq  GL(g,p) \rtimes {\rm Sym}_g (p)
$$
where Sym$_g (p)$ denotes the additive group of symmetric $g \times g$-matrices over $\F_p$.
This implies
\begin{eqnarray*}
\nu (p, \cdots ,p) &=& \frac{|Sp(2g,p)|}{|GL (n,p)| \cdot p^{\frac{g(g+1)}{2}}}\\[2ex]
&=& \frac{p^{g^2} \prod^g_{i=1} (p^{2i} -1)}{\prod^g_{i=0} (p^g - p^{i-1}) p^{\frac{g(g+1)}{2}}}\\[2ex]
&& \mbox{(by [H], II, Hilfssatz 6.2 and Satz 9.13).}\\[2ex]
&=& \prod^g_{i=1} (p^i +1).
\end{eqnarray*} \hfill $\square$
\medskip

\noindent
{\bf Proposition 4.3:} {\it $\nu(1,...,1,d) \,\, = \,\, \sigma(d) \,\,=\,\, \sum_{n|d} n$ for 
any positive integer $d$.}
\bigskip

\noindent
{\it Proof:} Note that $K(1,...,1,d) \,\,=\,\, K(d)$ and any subgroup of order $d$ of $K(d)$ 
is maximal isotropic. It is well-known that
$(\Z/d \Z)^2$ contains exactly $\sigma (d)$ subgroups of order $d$. \hfill $\square$\\

In the remaining cases it is a little more difficult to compute the number $\nu(d_1,...,d_g)$, mainly because 
it may happen that there are maximal isotropic subgroups of different types. We do this only for 
$K(d_1,d_2)$ which turns out to be sufficient in order to determine the number of equivalence
classes of curves of genus $3$ in any primitive cohomology class of an abelian threefold.\\

In order to compute $\nu(d,d)$ it suffices according to Proposition 4.1 
to compute $\nu(p^n,p^n)$ for every prime power $p^n$.
The different types of maximal isotropic subgroups of $K(p^n,p^n)$ are listed in the following
table together with a typical example. For the examples recall the standard generators of $K(p^n,p^n)$
with $e(f_1,f_3) \,\,=\,\, e(f_2,f_4) \,\,=\,\,exp({2 \pi i \over p^n})$,  
$e(f_3,f_1) \,\,=\,\, e(f_4,f_2) \,\,=\,\,exp(-{2 \pi i \over p^n})$ and $e(f_\nu, f_\mu) \,\,=\,\,1$ 
otherwise. In the following table (and only there) we denote for abbreviation $\Z_{p^m} := \Z/p^m\Z$ 
for all $m$.                   
  
\begin{center}
\begin{tabular}{c|c|c|c}
\quad type \quad & \quad $M$isomorphic to \quad & \quad restrictions \quad & \quad  example \\ \hline
$1$ & $\Z^2_{p^n}$ &  & $<f_1,f_2>$ \\ \hline
$2_k$ & $\Z_{p^n} \times \Z_{p^{n-k}} \times \Z_{p^k}$ & $0 < k < n-k< n$ & $<f_1,p^kf_2,p^{n-k}f_4>$\\ \hline
$3$ & $\Z_{p^n} \times \Z_{p^k} \times \Z_{p^k}$ &  $2k=n$ & $<f_1,p^kf_2,p^kf_4>$ \\ \hline
$4_{k,l}$ & $\Z_{p^{n-l}} \times \Z_{p^{n-k}} \times \Z_{p^k} \times \Z_{p^l}$ &  $0<l<k<n-k$ & 
                                                                      $<p^lf_1, p^kf_2, p^{n-k}f_4,p^{n-l}f_3>$\\ \hline
$5_l$ &  $\Z_{p^{n-l}} \times \Z_{p^k} \times \Z_{p^k} \times \Z_{p^l}$ & $0<l<k<n, 2k=n$ & 
                                                                      $<p^l f_1, p^k f_2, p^k f_4,p^{n-l} f_3>$ \\ \hline
$6_k$ & $\Z_{p^{n-k}} \times \Z_{p^{n-k}} \times \Z_{p^k} \times \Z_{p^k}$ & $0<k<n-k<n$ & 
                                                                      $<p^k f_1, p^k f_2, p^{n-k} f_3,p^{n-k} f_4>$ \\ \hline
$7$ & $\Z_{p^k} \times \Z_{p^k} \times \Z_{p^k} \times \Z_{p^k}$ & $ 2k=n$ & 
                                                                      $<p^k f_1, p^k f_2, p^k f_3,p^k f_4>$ \\ \hline                                                                     
\end{tabular}
\end{center}

\bigskip                           
\noindent
Note that there are some obvious restrictions: For example type $4_{k,l}$ only occurs for $n \geq 5$. For types $3$, $5_k$ and $7$ the
number $n$ is necessarily even.
 
The following proposition computes the number $\nu (p^n,p^n)(-)$ of maximal isotropic subgroups of $K(p^n,p^n)$ of type $-$ 
from the previous table.\\

\noindent
{\bf Proposition 4.4:} (1) \quad {\it $ \nu (p^n,p^n)(1) = p^{3n-3}(p^2 + 1)(p+1),\\
\smallskip
(2) \quad \nu (p^n,p^n)(2_k) = p^{3n-2k-4}(p^2 + 1)(p+1)^2,\\
\smallskip
(3) \quad \nu (p^n,p^n)(3)  = p^{2n-3}(p^2 + 1)(p+1),\\ 
\smallskip                                                                                           
(4) \quad \nu (p^n,p^n)(4_{k,l})  = p^{3n - 4l -2k - 4}(p^2+1)(p+1)^2,\\ 
\smallskip                          
(5) \quad \nu (p^n,p^n)(5_l) = p^{2n - 4l -3}(p^2 + 1)(p + 1),\\
\smallskip
(6) \quad \nu (p^n,p^n)(6_k) = p^{3n-6k-3}(p^2 + 1)(p+1)\\
\smallskip        
(7) \quad \nu (p^n,p^n)(7)  = 1$. } \\ 

\noindent
{\it Proof}: The number $\nu (p^n,p^n)(-)$ of maximal isotropic subgroups of $K(p^n,p^n)$ of type $-$ 
can always be computed as 
$$
\nu (p^n,p^n)(-) = {N(-) \over D(-)}
$$
where  $N(-)$ denotes the number of ordered bases of maximal isotropic subgroups of $K(p^n,p^n)$
of type $-$ and $D(-)$ denotes the number of ordered bases of the corresponding abelian
group.\\
There is a more elegant, though conceptually more involved proof of (1): 
Let $\tilde{e}_n : (\Z/p^n \Z)^4 \times (\Z/p^n \Z)^4 \lra \Z/p^n \Z$ 
denote the additive version of the multiplicative alternating form $e( \cdot , \cdot ): (\Z/p^n \Z)^4 \times 
(\Z/p^n \Z)^4 \lra \C^*$ of above.
The alternating forms $\tilde{e}_n$ form a projective system and define an alternating form $\tilde{e}: (\Z_p)^4 \times
(\Z_p)^4 \lra \Z_p$ (Here $\Z_p$ denotes the ring of p-adic integers). The symplectic group $Sp_4 (\Z_p)$ acts transitively 
on the set of maximal isotropic subgroups of $(\Z_p)^4$. Let $S$ denote the stabilizer of a fixed maximal isotropic subgroup. The 
set of all maximal isotropic subgroups of $(\Z_p)^4$ can be identified with the smooth scheme $M := Sp_4 (\Z_p)/S$ over $\Z_p$. 
The maximal 
isotropic subgroups of $(\Z/p^n \Z)^4$ can be considered as the $\Z/p^n \Z$-valued points of the scheme $M$. 
Since for any smooth scheme $N$ of relative dimension $d$ over $\Z^p$  
the $\Z/p^n \Z$-valued points and the $\Z/p \Z$-valued points of $N$ are related by
$$
\# N(\Z/p^n\Z) = p^{d(n-1)} \cdot \# N(\Z/p\Z)
$$
and since $M$ is of relative dimension 3 over $\Z_p$, Proposition 4.2 implies assertion (1).\\ 
Assertions (4), (5) and (6) can be deduced from (1), (2) and (3) using the following remark:
If $W \subset K(p^n,p^n)$ is any isotropic subgroup then the maximal isotropic subgroups containing $W$ are in natural bijection with 
the maximal isotropic subgroups of $W^\perp/W$. Thus we find: $\nu(p^n,p^n)(4_{k,l}) = \nu(p^{n-2l},p^{n-2l})(2_{k-l})$, 
$\nu(p^n,p^n)(5_l) = \nu(p^{n-2l},p^{n-2l})(3)$ and $\nu(p^n,p^n)(6_k) = \nu(p^{n-2k},p^{n-2k})(1)$.   \hfill $\square$

\bigskip 
\noindent
It is easy to make a similar computation for  $\nu (p^m,p^n)$ with $m<n$. However there are considerably 
more types of maximal isotropic subgroups to distinguish. We omit the corresponding tables and formulas. 
The following table gives the number $\nu (d_1,d_2)$ for small $d_1, d_2$, which are computed using
Propositions 4.1, 4.3, 4.4 and the corresponding formulas for $m<n$.                                                                                
        
\begin{center}
\begin{tabular}{c|c|c|c|c|c}
\quad $d$ \quad & \quad $\nu (1,d)$ \quad & \quad $\nu (d,d)$ \quad &  & \quad $(d_1,d_2)$ \quad & \quad $ \nu (d_1,d_2)$ \\ \hline
2 & 3 & 15 & & (2,4) & 51 \\ \hline
3 & 4 & 40 & & (2,6) & 60 \\ \hline
4 & 7 & 151 & & (2,8) & 114 \\ \hline
5 & 6 & 156 & & (2,10) & 90 \\ \hline
6 & 12 & 600 & & (2,12) & 204 \\ \hline
7 & 8 & 400 & & (3,6) & 120 \\ \hline
8 & 15 & 1335 & & (3,9) & 184 \\ \hline
9 & 13 & 1201 & & (3,12) & 280 \\ \hline
10 & 18 & 2340 & & (4,8) & 363 \\ \hline
11 & 12 & 1464 & & (4,12) & 604 \\ \hline
12 & 28 & 6040 & & (5,10) & 468 \\ \hline
16 & 31 & 10191 & & (6,12) & 2040 \\ \hline
\end{tabular}
\end{center}

\bigskip

\noindent
{\bf 5. Curves in a minimal cohomology class.}
\vspace{1ex}
\newline
\noindent
Let $(A,L)$ be a polarized abelian variety of type $(d_1,...,d_g)$. Let $c_1(L)$ denote the cohomology class
of the line bundle $L$. The class $\bigwedge^{g-1} c_1(L) \in H^{2g-2} (A,\Z)$ is divisible by $(g-1)!d_1 \cdots d_{g-1}$.
The class ${1 \over (g-1)!d_1 \cdots d_{g-1}} \bigwedge^{g-1} c_1(L)$ is not divisible in
$H^{2g-2} (A,\Z)$ and is called the {\it minimal cohomology class \,\,(of dimension $1$)} in $(A,L)$ (This follows easily
from the fact that for a suitable choice of a real basis $x_1, \ldots, x_{2g}$ of the tangent space of $A$ at 0 we have
$c_1(L) = - \sum\nolimits_{i=1}^g d_i dx_i \wedge dx_{g+i}$ (see [CAV], Lemma 3.6.4)). Denote by
$N_{min}(d_1,...,d_g)$ the number of translation classes of curves of genus $g$ in the minimal cohomology class
of $(A,L)$. If $L'$ is a $d_1$-th root of $L$, then $L'$ defines a polarization of type $(1,{d_2 \over d_1},...,{d_g \over d_1})$
and the minimal cohomology classes in $(A,L)$ and $(A,L')$ coincide. In particular we have $N_{min}(d_1,...,d_g) = 
N_{min}(1,{d_2 \over d_1},...,{d_g \over d_1})$. Hence we may always assume that $(A,L)$ is of type $(1,d_2,...,d_g)$.\\

In the case of an abelian surface $(A,L)$ of type $(1,d)$ the minimal cohomology class is the polarization $c_1(L)$ itself.
In this case the number $N_{min}(1,d)$ has been computed by G\"ottsche [G], Debarre [D] and Bryan - Leung [BL]. 
Theorem 2.1 and Proposition 4.3 yield
$$
N_{min}(d_1,d_2) = \sigma ({d_2 \over d_1})
$$
Note that the results of Section 4 also give the number of translation classes of curves of genus 2 in non minimal cohomology 
classes of abelian surfaces.\\

Now let $(A,L)$ be a simple abelian threefold of type $(1,d_2,d_3)$. Recall from Section 4 that $\nu (m,n)$ denotes the number of 
maximal isotropic subgroups of $K(m,n)$, hence it equals the number of maximal isotropic subgroups of 
$K(1,m,n)$. So Corollary 2.3 yields
$$
N_{min}(1,d_2,d_3) = \nu({d_3 \over d_2}, d_3)     \eqno(4)
$$
and we obtain from Proposition 4.3:\\

\noindent
{\bf Proposition 5.1:} {\it For any simple abelian threefold $(A,L)$ of type $(1,d,d)$ we have
$$
N_{min}(1,d,d) =  \sigma (d).
$$ }

\noindent
Similarly we have\\       

\noindent                                                                                             
{\bf Proposition 5.2:} {\it Let $(A,L)$ be a simple abelian threefold of type $(1,1,d)$ and $d \, = \, \prod^r_{i=1} p_i^{n_i}$
the prime decomposition. Then
$$ 
N_{min}(1,1,d) = \prod^r_{i=1} \nu (p^{n_i}_i, p^{n_i}_i)
$$
and for any prime number $p$ we have 
\begin{eqnarray*}
\nu (p^{2m+1},p^{2m+1}) &=& (p^2 + 1)(p+1)[p^{6m} + {p^{4m-1}(p^{2m}-1) \over p-1} + {p^{6m} - 1 \over p^6 - 1}] +\\
&& + {p+1 \over p-1}[ p^{4m-1}{(p^{2m-2} - 1) \over p^2 - 1} - p^3 {p^{6m-6} - 1 \over p^6 - 1}]
\end{eqnarray*} 
for all $m \geq 2$ and
\begin{eqnarray*}
\nu (p^{2m},p^{2m}) &=& (p^2 + 1)(p+1)[p^{6m-3} + p^{4m-2} {p^{2m-2} - 1 \over p-1} + p^{4m-3} + p {p^{4m-4}-1 \over p^4 -1} + p^3 {p^{6m-6} -1
\over p^6 -1}] \\
&& + {p+1 \over p-1} [ p^{4m-2} {p^{2m-4} - 1 \over p^2 - 1} - p^6 {p^{6m-12} - 1 \over p^6 - 1} ] + 1                                         
\end{eqnarray*}
for all $m \geq 3$. }\\

\noindent
For smaller values of $m$ the formulas are similar, but simpler.\\

\noindent
{\it Proof}: By identity (4) we have $N_{min}(1,1,d) = \nu (d,d)$ and by applying Proposition 4.1 we obtain 
$ \nu (d,d) = \prod_{i=1}^r \nu (p_i^{n_i},p_i^{n_i})$. Thus we get
$$
N_{min}(1,1,d) = \prod_{i=1}^r \nu(p_i^{n_i},p_i^{n_i}).
$$
The expression for $\nu(p^{2m-1},p^{2m-1})$ and $\nu (p^{2m},p^{2m})$ is obtained by adding the appropriate terms in 
Proposition 4.4. \hfill $\square$\\

Similar formulas can be given for $N_{min}(1,d_2,d_3)$ for $d_2 \,\,<\,\, d_3$. We omit them since they look more complicated. 
But note that the table on page 9 gives $N_{min}(1,d_2,d_3)$ in some cases. For example $N_{min}(1,2,4)\,\,=\,\, 51$ etc.

\vspace{1cm}
\begin{center}
{\bf References:}
\end{center}
\begin{itemize}
\item[{[BL]}] Ch. Birkenhake, H. Lange: The dual polarization of an abelian variety, Archiv Math. 73 (1999), 380 - 389.
\item[{[BrL]}] J. Bryan, N.C. Leung: Generating functions for the number of curves on abelian surfaces, Duke Math. 
J. 99 (1999), 311-28.
\item[{[D]}] O. Debarre: On the Euler characteristic of generalized Kummer varieties, Am. J. Math. 121 (1999), 577-86. 
\item[{[G]}]  L. G\"{o}ttsche: A conjectural generating function for numbers of curves on surfaces,
Comm. Math. Physics 196 (1998), 523-33.
\item[{[H]}] B. Huppert: Endliche Gruppen, Band 1, Springer, 1967.
\item[{[CAV]}] H. Lange, Ch. Birkenhake: Complex Abelian Varieties, Grundlehren Nr. 302, Springer Verlag, (1992)
\item[{[M]}] T. Matsusaka: On a characterization of a Jacobian Variety. Mem. Coll. Science, Kyoto,
Ser. A. 32 (1959), 1 - 19.
\end{itemize}

\vskip1cm\noindent
Mathematisches Institut
\hfill\break
Bismarckstr. $1{1\over 2}$, D-91054 Erlangen, Germany\\
lange@mi.uni-erlangen.de

\bigskip\noindent
Dipartimento di Matematica, Universit\`a Roma III
\hfill\break
L.go S.L. Murialdo 1, 00146 Roma (Italy)\\
sernesi@mat.uniroma3.it
\end{document}